\documentclass{amsart}
\usepackage{amssymb}
\usepackage{tikz}
\usepackage{hyperref}

\usepackage{epic,curves,mfpic}
\newtheorem{theo}{Theorem}[section]

\newtheorem{conj}[theo]{Conjecture}

\newtheorem{prop}[theo]{Proposition}

\theoremstyle{remark}
\newtheorem{exam}[theo]{Example}

\newcommand{\Z}{\mathbb{Z}}
\newcommand{\R}{\mathbb{R}}

\newcommand{\Y}{\mathcal{Y}}

\newcounter{fig}
\begin{document}

\title{ Finitely Presented Groups Acting on Trees}

\author{M.J.Dunwoody}

\subjclass[2010]{20F65 ( 20E08)}
\keywords{groups acting on trees, group splittings}

\begin{abstract}
It is shown that for any action of a finitely presented group $G$ on an $\R$-tree, there is a decomposition of $G$ as the fundamental
group of a graph of groups related to this action.    If the action of $G$ on $T$ is non-trivial, i.e. there is no global fixed point, then $G$ has
a non-trivial action on a simplicial $\R $-tree.
\end{abstract}
\maketitle

\begin{section} {Introduction}

 A group $G$ is said to split over a subgroup $C$ if either $G = A*_C B$, where $A \not= C$
 and $B \not= C$ or $G$ is an HNN-group $G = < A *_C = < A, t | t^{-1}at = \theta (a)>$ where 
 $\theta : C \rightarrow A$ is an injective homomorphism.
 It is one of the basic results of Bass-Serre theory (see \cite{[DD]} or [\cite{[Se]}), that a finitely generated group $G$ splits
 over some subgroup $C$ if and only if there is an action of $G$ on a  tree $T$, without inversions,
 such that
 for no vertex $v \in VT$ is $v$ fixed by all of $G$.   Here the tree is a combinatorial tree, i.e. a connected 
 graph with no cycles, and an action without inversions is one in which no element $g \in G$ transposes the
 vertices of an edge.  
 Tits \cite{[T]} introduced the idea of an $\R $-tree, which is a non-empty metric space in which any two points
 are joined by a unique arc, and in which every arc is isometric to a closed interval in the real line $\R $.
 Alternatively an $\R $-tree is a $0$-hyperbolic space.   A tree in the combinatorial sense
 can be regarded as a $1$-dimensional simplicial complex.   The polyhedron of this complex will be an
 $\R $-tree - called a simplicial $\R $-tree.  However not every $\R $-tree is like this.
 A point $p$ of an $\R $-tree $T$ is called {\it regular} if $T - p$ has two components.
 An $\R $-tree is simplicial if the points of $T$ which are not regular form a discrete subspace of 
 $T$.  It is fairly easy to construct examples of $\R $-trees where the set of non-regular points is not
 discrete.
 There are good introductory accounts of groups acting on $\R $-trees in \cite{[Be]} and \cite{[Sh]}.
 We assume that all our actions are by isometries.
 It is a classical result that a group is free if and only if it has a free action on a simplicial tree.
 As the real line $\R $ is an $\R $-tree and $\R $ acts on itself freely by translations, any free abelian group has a free action on a $\R $.
 Morgan and Shalen \cite{[MS]} showed that the fundamental group of any compact surface other than
 the projective plane and the Klein bottle has a free action on an $\R $-tree. Rips  showed that the only finitely generated groups that act freely on an $\R $-tree are 
 free products of free abelian groups and surface groups.  Rips never published his proof, but there are
 proofs of more general results by Bestvina -Feighn \cite{[BF]} and by Gaboriau-Levitt-Paulin (see \cite{[P]} or \cite{[C]}).
 Bestvina and Feighn classify the {\it stable } actions of finitely generated groups on $\R$-trees. 
 Recall, that an action of a group $G$  on an $\R $-tree    is said to be {\it stable} if there is no sequence of arcs $l_i$ such that $l_{i+1}$
is properly contained in $l_i$ for every  $i$, and for which the stabilizer $G  _I$ of $l_i$ is properly contained in $G  _{i+1}$ for every $i$.
In particular \cite{[BF]} Bestvina and Feighn proved that if a finitely presented group has a non-trivial minimal stable action on an $\R$-tree then it has a non-trivial 
action on some simplicial tree.

 A group is said to be $(FA)$ if it has no non-trivial action on a simplicial $\R$-tree and it is said to be $(F\R)$ if it has no non-trivial action on
any $\R $-tree.  A {\it trivial} action is one in which there is a point of the tree that is a global fixed  point.
   In contrast  A.Minasyan \cite {M} and I ~\cite{[D2]} in separate papers have  given examples of finitely generated groups that are
$(FA)$ but not $(F\R )$.  These provided a negative answer to Shalen's Question A of \cite{[Sh]}.     In an earlier paper \cite {[D2]} I gave an example of a finitely generated group
that had a non-trivial action on an $\R$-tree with finite cyclc arc stabilizers but for which any simplicial decomposition has an edge group that contains a non-cyclic free group.
This gave a negative answer to Conjecture D of \cite {[Sh]}.
In this paper it is shown that there are positive  answers to these questions  for finitely presented groups.   The situation is therefore similar to that of accessibility in finitely generated groups, in that finitely presented groups are accessible \cite {[D1]}, but there are examples of finitely generated groups that are not accessible \cite {[D3]},\cite {[D4]}.   The questions are closely related.

A {\it morphism}  from a segment $I$ to an $\R$-tree $T$ is a continuous map $\phi  : I \rightarrow T$ such that $I$ may be subdivided
into finitely many subsegments that $\phi $ maps isometrically into $T$.   Let $T, T'$ be $\R $-trees with actions of groups $G, G'$ respectively.  Let
$\rho : G\rightarrow G'$ be a homomorphism.   A {\it morphism } from $T$ to $T'$ is a map $\phi $ equivariant with respect to $\rho $ which induces
a morphism on every segment $I \subset T$.

In this paper  the following theorem is proved.

\begin{theo}\label {main}

 Let $G$ be a finitely presented group
and let $T$ be a $G$-tree, i.e. an $\R$-tree on which $G$ acts by isometries.

Then 
$G$ is the fundamental group of a finite graph $(\Y, Y)$  of groups, in which every edge group  is finitely generated and fixes a point of $T$.   If $v \in VY$, then either $\Y (v)$ fixes a vertex of $T$ or there is a homomorphism from $\Y (v)$ to a target group  $Z(v)$ (a parallelepiped group), which is the fundamental group of a cube complex of groups based on a single $n$-cube $c(v)$.

 Every hyperplane of $c(v)$ is associated with a non-trivial splitting of $G$.

There is a marking of  the  cube $c(v)$ so that the corresponding $\R $-tree with its $Z(v)$-action
is the image of a morphism from a $\Y (v)$-tree $T_v$ and this tree is the minimal $\Y (v)$-subtree of $T$.
\end{theo}

The action of a target group on an $\R $-tree is usually unstable, but a parallelepiped group of rank $n$ contains a free abelian group of rank $n$
and this acts freely on $\R $ by translation. 

The main theorem in a previous version of this paper is incorrect.  I thought that the Higman group $H = \langle  a, b, c, d | aba^{-1} = b^2, bcb^{-1} = c^2, cdc^{-1} = d^2, dad^{-1} = a^2\rangle $.
 had an action on a nonsimplicial $R$-tree, since it had non-compatible decompositions as a free product with amalgamation.
 In fact this is not the case.   I thought that adding linear combinations of tracks corresponding to these decompositons would result in infinitely many non-trivial such
 decompositions.   In fact all the tracks obtained correspond to trivial decompositions (see \cite {BD}).
 
 We do not encounter Levitt (or thin) type actions in our analysis (see \cite {[Be]}). This is because any such action is resolved by a simplicial action.


\end{section}
\begin{section}{Target Groups}

In \cite{[D5]}  {\it rectangle groups}  were constructed.  

The rectangle group $R = R(m.n.p.q),  m,n,p,q \in \{2,3, \dots \}\cup \{\infty \}$
is the group with presentation 
$$R = \{ a, b, c, d | a^m = b^n = c^p = d^q =1, ab^{-1} = cd^{-1}, ac^{-1} = bd^{-1} \}. $$

  \begin{figure}[htbp]
\centering
\begin{tikzpicture}[scale=.7]

  \draw  (0,.0) -- (4, 0) -- (4, 3) -- (0, 3)-- (0,0)  ;
{[blue]
\draw [left] (0,0) node {$a$} ;
\draw [right] (4,0) node {$c$} ;
\draw [left] (0,3) node {$b$} ;
\draw [right] (4,3) node {$d$} ;
}
  \draw  (6,.0) -- (10, 0) -- (10, 3) -- (6, 3)-- (6,0)  ;
{[red]
\draw [left] (6,0) node {$m$} ;
\draw [right] (10,0) node {$p$} ;
\draw [left] (6,3) node {$n$} ;
\draw [right] (10,3) node {$q$} ;

}
\end{tikzpicture}
\vskip-3mm\caption{rectangle group }\label{fig:mark}

\end{figure}
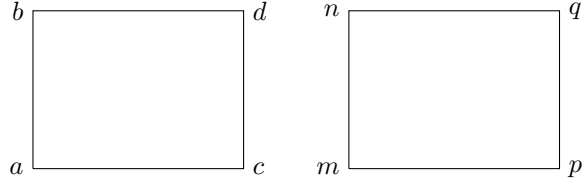

Think of the relations as saying that opposite edge vectors are equal and that the corners are assigned 
orders,  a corner can have infinite order. 

In the group $R$ above, let $x = ab^{-1} = cd^{-1},  y = ac^{-1} = bd^{-1}$,
then $xy = ab^{-1}bd^{-1} =ad^{-1} = ac^{-1}cd^{-1} = yx$,  and $x, y$ generate a free abelian rank $2$ group.

Also $R$ has {\it incompatible } decompositions as a free product with amalgamation
$$R = \langle  a, b \rangle * _{ab^{-1}  = cd^{-1}, b^{-1}a = d^{-1}c} \langle c ,d \rangle ,$$
and
 $$R = \langle  a, c \rangle * _{ac^{-1}  = bd^{-1}, c^{-1}a = d^{-1}b} \langle b ,d \rangle ,$$

the amalgamated subgroup in each case is free of rank two.

The decompositions are incompatible because the lines dividing the rectangle intersect.
 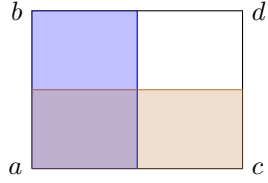
\begin{figure}[htbp]
\centering
\begin{tikzpicture}[scale=.7]

  \draw  (0,.0) -- (4, 0) -- (4, 3) -- (0, 3)-- (0,0)  ;
{[red]
\draw [left] (0,0) node {$a$} ;
\draw [right] (4,0) node {$c$} ;
\draw [left] (0,3) node {$b$} ;
\draw [right] (4,3) node {$d$} ;
}
\draw [blue]  (2,0)--(2,3) ;
\draw [brown] (0, 1.5) --(4, 1.5) ;
\filldraw [blue] [opacity = .25] (2,0)--(2,3)--(0,3)--(0,0)--cycle ;
\filldraw [brown] [opacity = .25] (0,1.5)--(4,1.5)--(4,0)--(0,0)--cycle ;

\end{tikzpicture}
\vskip-3mm\caption{Decompositions of  a rectangle group\label {rect}}

\end{figure}

A {\it cube complex} is similar to a simplicial complex except that the building blocks are $n$-cubes rather
than $n$-simplexes.    A rectangle group $R$ acts on a simply connected $2$-dimensional cube complex $\tilde C$ with orbit space $C$.
This is illustrated in Fig~\ref{lattice}.
Apart from some exceptional cases, when two or more of $m, n, p, q$  are $2$, there are three orbits of $2$-cells, each with trivial stabilizer.  In all cases  there are  four orbits of edges also with trivial stabilizers, and four orbits
of vertices labelled $A, B, C, D$ with stabilizers which are cyclic of orders $m, n, p, q$ respectively.
In the group $R(2,2,2,2)$ the subgroup $\langle x, y \rangle$ has index $2$ and there is one orbit of $2$-cells.
In $R(2,2,2,q)$ or $R(2,2,p, q)$ for $p, q \geq 3$ there are two orbits of $2$-cells and both $C$ and $\tilde C$ are  $2$-orbifolds.

In \cite{[D5]}  it is shown that for any action of $J = \langle x, y\rangle$ on an $\R $-tree there is an action of  the rectangle group on an $\R $ tree $T$ which restricts to the given action on the minimal $J$-subtree of $T$.   This action is unstable in all cases when there are $3$ orbits of $2$-cells.

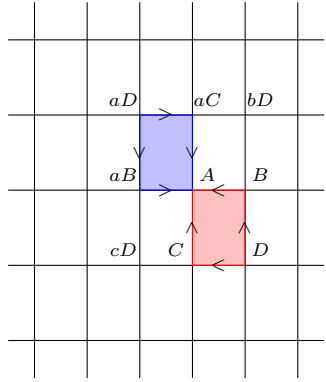
\begin{figure}[htbp]

\begin{tikzpicture}[xscale = .07, yscale = .1, fill opacity =.25]
\centering

\draw (0,0)--(0,50) ;
\draw  (10, 0 )--  (10,50) ;
\draw  (20, 0 )--  (20,50) ;
\draw  (30, 0 )--  (30,50) ;
\draw  (40, 0 )--  (40,50) ;
\draw  (50, 0 )--  (50,50) ;

\draw (-5,5)--(55,5) ;
\draw (-5,15)--(55,15) ;
\draw (-5,25)--(55,25) ;
\draw (-5,35)--(55,35) ;
\draw (-5,45)--(55,45) ;

\filldraw [red]  (30, 15)--(40, 15)--( 40, 25) --(30,25) -- cycle  ;
\filldraw [blue]   (20,25)--(30, 25)--( 30, 35) --(20,35) -- cycle  ;

\draw [fill opacity =1] [left, above] (17, 15) node  {$_{cD}$} ;
\draw [fill opacity =1] [left,above] (17, 25) node  {$_{aB}$} ;
\draw [fill opacity =1][left, above] (17, 35) node  {$_{aD}$} ;
\draw [fill opacity =1, left, above ](27, 15) node  {$_C$} ;
\draw [fill opacity =1,left,above ](43, 15) node  {$_D$} ;
\draw [fill opacity =1,left,above](43, 25) node  {$_B$} ;
\draw [fill opacity =1,left,above ](33, 25) node  {$_A$} ;
\draw [fill opacity =1,left,above](33, 35) node  {$_{aC}$} ;
\draw [fill opacity =1,right,above](43, 35) node  {$_{bD}$} ;

\draw [fill opacity =1](25, 35) node  {$_>$} ;
\draw [fill opacity =1](35,15) node  {$_<$} ;
\draw [fill opacity =1](25, 25) node  {$_>$} ;
\draw [fill opacity =1](35,25) node  {$_<$} ;
\draw [fill opacity =1](30, 30) node  {$_\vee$} ;
\draw [fill opacity =1](20,30) node  {$_\vee$} ;

\draw [fill opacity =1](30, 20) node  {$_\wedge$} ;
\draw [fill opacity =1](40,20) node  {$_\wedge$} ;


\end{tikzpicture}
\caption{The Euclidean space for a rectangle group}\label{lattice}
\end{figure}

There is a Euclidean $2$-dimensional subspace $E$  of $\tilde C$ acted on by $\langle x, y \rangle$ .
For the action of $\langle x, y \rangle $ on $E$ there is one orbit of $2$-cells, each of which is made up of $4$ smaller rectangles
of $\tilde C$.  In the diagram the points $A, B, C, D$ are stabilised by $a, b, c, d$ respectively.

Note that the blue rectangle $A, aB, aC, aD$ is in the same $R$-orbit as the red rectangle $A, B, C,D$, and $bD = bd^{-1}D =  yD,
cD = cd^{-1}D =xD$.

A {\it parallelepiped group} of dimension $n$ has $2^n$ generators corresponding to the vertices of an $n$-cube.   The generators
corresponding to a $2$-dimensional face satisfy the relations of a rectangle group.
Such a group has an action on an $n$-dimensional cube complex $C_n$ for which the orbit space is an $n$-cube.
There is a subgroup $J_n$ that is free abelian of rank $n$, which acts on a subcomplex  $E_n$ of $C_n$ so that the orbit
space $J_n\backslash E_n$ consists of $2^n$ smaller cubes. 

In \cite {[D5]} there is a detailed description of the action for $n=3$.

\end{section}
\

\begin{section}{ Finitely presented groups}

Let $X$ be a  finite CW $2$-complex.
We introduce the idea of a complex of
groups $G(X)$ based on
$X$.  This is a slightly different notion to a
special case of the complex of groups described
by
 Haefliger \cite  {[Ha]}.  Haefliger restricts
$X$ to be a simplicial cell complex.
One can get from our situation to that of
Haefliger by triangulating each $2$ cell.
We are only concerned with the situation
when each group assigned to a $2$-cell  is
trivial.

Thus the $1$-skeleton $X^{1}$ of $X$ is
a graph. We take the edges to be oriented,
and use Serre's notation, so that each
edge $e$ has an initial vertex $\iota e$
and a terminal vertex $\tau e$ and $\bar
e$ is $e$ with the opposite orientation.
Let $G(X^1)$ be a graph of groups based
on $X^1$
The attaching map of each $2$-cell $\sigma $
is given by a closed path in $X^1$.
Let $S$ be a spanning tree in $X^1$.
The fundamental group $\pi (G(X), S)$ of
the complex of groups $G(X)$ is the fundamental
group of the graph of groups $G(X^1)$
together with  extra relations corresponding
to the attaching maps of the $2$-cells.
Thus $\pi (G(X), S)$ is generated
by the groups $G(v), v \in V(X^1)$ and
the elements $e \in E(X^1)$.
For each $e \in E(X^1),
G(e)$ is a distinguished subgroup of $G(\iota
e)$ and there are injective homomorphisms $t_e
: G(e)
\rightarrow G(\tau e), g \mapsto g^{\tau e}$.
The relations of $\pi (G(X), S)$ are
as follows:-

the relations for $G(v),$ for each $v \in
V(X^1)$

$e^{-1}ge = g^{\tau e}$ for all $e \in E(X^1),
g \in G(e) \leq G(\iota e),$ 

$e =1$ 
if $e \in
E(S)$.

For each attaching closed path $e_1, e_2, \dots
, e_n$ in $X$ of a $2$-cell, there is a relation
$$g_0e_1g_1e_2x_2\dots g_{n-1}e_n = 1,$$ where
$g_i \in G_{\tau e_i} = G_{\iota e_{i+1}}$, called the {\it attaching word}.  The elements $g_i$ are called {\it joining } elements.
Such a word represents both a path $p$, called the {\it attaching path} in the  Bass-Serre tree $T$ corresponding to 
the graph of groups $G(X^1)$, for which initial  point $\iota p$ and end point $\tau p$ are in the same
$\pi (G(X^1),S)$-orbit and an element $g \in \pi(G(X^1), S)$ for which $g\iota p = \tau p$.
Adding the relation identifies the points $\iota p$ and $\tau p$ and puts $g = 1$.
If we carry out all these identifications, we obtain a $G$-graph $\Gamma $
in which the attaching paths are all closed paths.
We describe specifically how this path arises (as in \cite {[DD]}, p15).
We lift $S$ to an isomorphic subtree of $S_1$ of $\Gamma $.  Thus the vertex set of $S_1$ is a transversal
for the action of $G$ on $\Gamma $.   For each edge $e$ in $X - S$ we can choose an edge 
$\tilde e \in T$ such that $\tilde e$ maps to $e$ in the natural projection, and $\iota \tilde e$
is a vertex of $S_1$.   Let $\tilde S  $ be the union of $S_1$ with these extra edges.  Note that it will not normally be the case that $\tau \tilde e \in \tilde S$ and so $\tilde S$ is not usually a subtree of $T$, but there will
be an element $c(e) \in G$ such that $c(e)^{-1}(\tau \tilde e) \in \tilde S$.  These elements (called the {\it connecting} elements) together with
the stabilizers of elements of $VS_1$ generate $G$.     
Clearly $\tilde S$ consists of a transversal for the action of $G$ on both the edges and vertices of $\Gamma $.
Let $\iota p = v_0$ be the vertex of $\tilde S$ lying above $\iota e_1$, and put $x_0 = g_0$   Suppose we have constructed
$v_i$ and $x_i \in G$ so that $v_i$ is  the terminal vertex of the path corresponding to $g_0e_1g_1e_2g_2\dots g_{i-2}e_{i-1}$, and
so that if   $\tilde v_i$ is  the element of $\tilde S$ in the orbit of $v_i$,  then $v_i = x_i\tilde v_i$.
This is certainly true when $i = 0$.
To construct $v_{i+1}$ and $x_{i+1}$, put $x_{i+1} = x_ic(e_{i+1})g_{i+1}$ where we put $c(e) = 1$ if $e \in S$.   Then $x_{i+1}\tilde v_{i+1} = x_ic(e_{i+1})\tilde v_{i+1}$ is the terminal vertex of the edge
$x_i\tilde e_{i+1}$ with initial vertex $v_i$.  Note that $x_i$ is obtained from $g_0e_1g_1e_2g_2\dots g_{i-2}e_{i-1}$ by replacing each $e_i$ by $c(e_i)$.

We now foliate each $2$-cell of $X$ in a
particular way.  Thus let 
$ D = \{ (x, y) | x, y \in {\bf R}, x^2 + y^2
\leq 1\}$
be the unit disc.

\vfill
\vskip-3mm
\begin{figure}[htbp]
\centering

\begin{tikzpicture}
          
    \path (0,0) coordinate (p1);
    \path (2, 2) coordinate (p2);
    \path (4,3) coordinate (p3);
    \path (6,2) coordinate (p4);
    \path (8,0) coordinate (p5);
    \path (6,-2) coordinate (q5);
    \path (2, -2) coordinate (q2);
    \path (4,-3) coordinate (q3);
    \path (6,-2) coordinate (q4);

     \draw (p1) -- (p2) -- (p3) --(p4) -- (p5)--(q4)--(q3)--(q2)--(p1) ;

    \draw (7.8, .2) --(7.8,  -.2) ;
\draw (7.6,.4)--(7.6,-.4) ;
\draw (7.4, .6) --(7.4,  -.6) ;
\draw (7.2,.8)--(7.2,-.8) ;
\draw (7, 1) --(7,-1) ;

\draw (6.8,1.2)--(6.8,-1.2) ;
\draw (6.6,1.4)--(6.6,-1.4) ;

\draw (6.4, 1.6) --(6.4, -1.6) ;
\draw (6.2, 1.8)--(6.2,-1.8) ;
\draw (6, 2) --(6,-2) ;

\draw (5.8,2.1)--(5.8,-2.1) ;
\draw (5.6,2.2) --(5.6,-2.2) ;
\draw (5.4,2.3)--(5.4,-2.3) ;
\draw (5.2,2.4) --(5.2,-2.4) ;

\draw (5,2.5)--(5,-2.5) ;
\draw (4.8,2.6) --(4.8,-2.6) ;
\draw (4.6,2.7)--(4.6,-2.7) ;
\draw (4.4,2.8) --(4.4,-2.8) ;
\draw (4.2,2.9) --(4.2,-2.9) ;
\draw (4, 3)--(4,-3) ;
   
    \draw (.2, .2) --(.2,  -.2) ;
\draw (.4,.4)--(.4,-.4) ;
\draw (.6, .6) --(.6,  -.6) ;
\draw (.8,.8)--(.8,-.8) ;
\draw (1, 1) --(1,-1) ;

\draw (1.2,1.2)--(1.2,-1.2) ;
\draw (1.4,1.4)--(1.4,-1.4) ;

\draw (1.6, 1.6) --(1.6, -1.6) ;
\draw (1.8, 1.8)--(1.8,-1.8) ;
\draw (2, 2) --(2,-2) ;

\draw (2.2,2.1)--(2.2,-2.1) ;
\draw (2.4,2.2) --(2.4,-2.2) ;
\draw (2.6,2.3)--(2.6,-2.3) ;
\draw (2.8,2.4) --(2.8,-2.4) ;

\draw (3,2.5)--(3,-2.5) ;
\draw (3.2,2.6) --(3.2,-2.6) ;
\draw (3.4,2.7)--(3.4,-2.7) ;
\draw (3.6,2.8) --(3.6,-2.8) ;
\draw (3.8,2.9) --(3.8,-2.9) ;

\draw (.8,.9)--(1, 1) ;
\draw (.9,.8)--(1,1);
\draw (2.8,2.5)--(3,2.5)--(2.9,2.35) ;
\draw (4.8,2.5)--(5,2.5)--(4.9, 2.7) ;
\draw (6.8, 1.1)--(7,1)--(6.9,1.2) ;
\draw (.8,-.9)--(1,-1) ;
\draw (.9,-.8)--(1,-1) ;
\draw (2.8,-2.5) --(3, -2.5)--(2.9,-2.35) ;
\draw (4.8, -2.5)--(5,-2.5)--(4.9,-2.7) ;
\draw (6.8, -1.1)--(7, -1)--(6.9, -1.2) ;

\draw  (1, 1.3) node {$e_1$} ;
\draw (1.9, 2.3) node  {$x$} ;
\draw (-.15, 0) node  {$v$} ;
\draw (1.9, -2.3) node{$y$} ;
\draw (1,-1.3) node {$e_2$} ; 

  \filldraw (p2) circle (1pt);
  \filldraw (p3) circle (1pt);  
 
  \filldraw (p4) circle (1pt);
    \filldraw (p5) circle (1pt) ;

  \filldraw (q2) circle (1pt);
  \filldraw (q3) circle (1pt);  
 
  \filldraw (q4) circle (1pt);
    \filldraw (p1) circle (1pt) ;

\end{tikzpicture}
\vskip-3mm\caption{Foliated $2$-cell}\label{fig:2cell}
\end{figure}
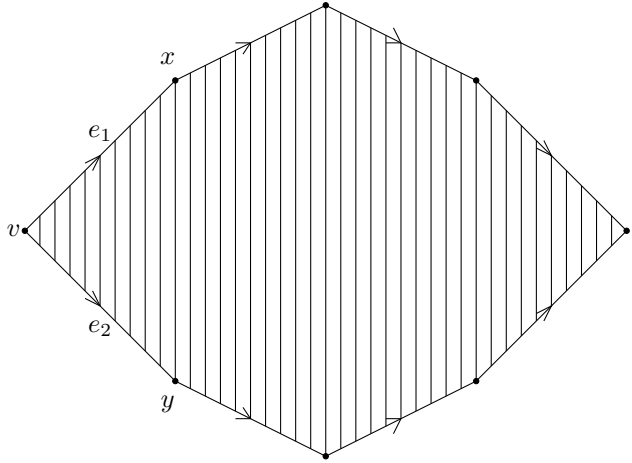

Give this the foliation in which leaves are
the intersection of $D$ with the vertical lines
$x = c$ where c is a constant in the interval
$[-1, 1]$.  
Let $\sigma $ be a $2$-cell of $X$ which is
attached via the closed path
$e_1, e_2, \dots , e_n$.  We map
$D$ to $\sigma $ so that for some $j = 2, \dots
, n-1$ the upper semi-circle joining
$(-1,0)$ and $(1,0)$ is mapped to the
path $e_1, \dots e_j$.
Thus there are points $z_0 = (-1, 0),
z_1, \dots , z_j = (1, 0)$ on the upper
semi-circle so that 
$z_i \mapsto \iota e_i, i = 1, 2, \dots
e_{j+1}$ and the map is continuous and
injective on each segment $[z_i, z_{i+1}]$,
except if $\iota e_i = \iota e_{i+1}$ in which
case the map is injective on the interior
points of this segment.
In a similar way the lower semi-circle
is mapped to the path $\bar e_n, \dots
,\bar e_{j+1}$.

Let, then, $X$ be a $2$-complex of groups
in which each $2$-cell is foliated as described
above and let $T$ be a $G$-tree, i.e. $T$
 is an ${\bf R}$-tree on which
$G$ acts by isometries.
We say that the $X$ {\it resolves} $T$
if there is an isomorphism 
$\theta : \pi (X, S) \rightarrow G$
which is injective on vertex groups
(and hence on all groups $G_{\sigma }$ for
all cells $\sigma $ of $X$).
In this situation (see \cite {[Ha]}),
the complex of groups is developable,
i.e. there is a cell complex $\tilde X$
on which $G$ acts and $G(X)$ is the complex
of groups associated with this action.
We also require that there be
a $G$-map  $\alpha : \tilde X \rightarrow T$
such that for each $1$-cell $\gamma $ the
restriction of $\alpha $ to $\gamma $ is
injective and for each $2$-cell
$\sigma $ and each $t\in T$, the intersection
of $\sigma $ with $\alpha ^{-1}(t)$ is either
empty or a leaf of the foliation described
above.

We show that if $G$ is finitely presented
then any $G$-tree has a resolution, i.e.
there is a cell complex $X$ as above that
resolves $T$.  Our approach is similar to that of \cite {[LP]}.

Since $G$ is finitely presented,
there is
simplicial $2$-complex
$X$ such that $\pi (X, S) \cong G$. Here
$S$ is a spanning tree in the
$1$-skeleton of $X$.
Let $\tilde X$ be the universal cover
of $X$.  Clearly there is a $G$-map
$\theta _0: V\tilde X \rightarrow T$,
which can be obtained by first mapping
a representive of each $G$-orbit of vertices
into $T$ and then extending so as to make
the map commute with the $G$-action.
Now extend this map to the $1$-skeleton
so that each $1$-simplex $\gamma $
with vertices $u, v$
of $\tilde X$ is mapped injectively to the
geodesic joining $\theta _0(u)$ and
$\theta _0(v)$. 
It may be necessary to subdivide $X$ and choose
the map $\theta _0$ to ensure that $\theta
_0(u) \not= \theta _0(v)$ for every
$1$-simplex $\gamma $.  We can extend the map
to every $1$-simplex so that it
commutes with the
$G$-action giving a $G$-map $\theta _1 :\tilde
X^1 \rightarrow T$.
Now we extend the map to the $2$-simplices.
Let $\sigma $ be a $2$-simplex with vertices
$u, v, w$.  If $\theta _0(u)$ lies on
the geodesic joining $\theta _0(v)$
and $\theta _0(w)$ then we can map
$\sigma $ as indicated in Fig~{\ref {simp}} (ii).
Each vertical line is mapped to a point.
If $\theta _0(u), \theta _0(v)$ and 
$\theta _0(w)$ are situated as in 
Fig~{\ref {simp}}(i) so that no point is on the geodesic
joining the other two, then we subdivide 
$\sigma $ as in Fig~{\ref {simp}} (iii).  The new vertex 
is mapped to the point $p$ of (i) and
the three new simplexes now have the middle 
vertex mapped into the geodesic joining the
images of the other two sides and are mapped
as shown in (iii).

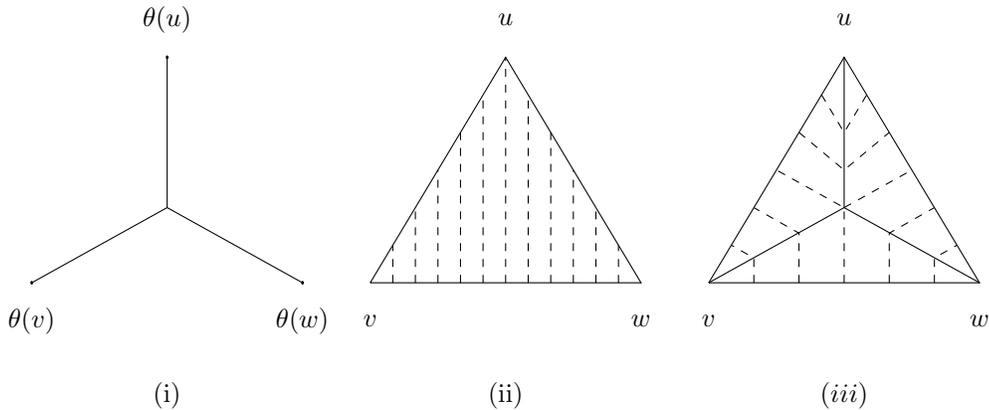
\begin{figure}

\centering
\begin{tikzpicture}[xscale=.3, yscale = .5]

     \filldraw (0,0) circle (1pt);
  \filldraw (6,6) circle (1pt);  
 
  \filldraw   (12,0) circle (1pt);
\draw  (15,-1)  node {$v$} ;
\draw  (21,7)  node {$u$} ;
\draw  (27,-1)  node {$w$} ;
\draw  (30,-1)  node {$v$} ;
\draw  (36,7)  node {$u$} ;
\draw  (42,-1)  node {$w$} ;
\draw  (6,7)  node {$\theta (u)$} ;
\draw  (0,-1)  node {$\theta (v)$} ;
\draw  (12,-1)  node {$\theta (w)$} ;
\draw  (6,-3)  node {(i)} ;
\draw  (21,-3)  node {(ii)} ;
\draw  (36,-3)  node {$(iii)$} ;



\draw  (0,0)--  (6,2) --  (6,6)   ;
\draw  (6,2)--(12,0)  ;
\draw  (15,0) --  (27,0) --  (21,6) -- (15,0) ;
\draw (30,0) --  (42,0)-- (36,6)--  (30,0) ;

\draw [dashed]  (16,0) --(16,1) ;
\draw [dashed]  (17,0) --(17,2) ;
\draw [dashed]  (18,0) --(18,3) ;
\draw [dashed]  (19,0) --(19,4) ;
\draw [dashed]  (20,0) --(20,5) ;
\draw [dashed]  (21,0) --(21,6) ;
\draw [dashed]  (22,0) --(22,5) ;
\draw [dashed]  (23,0) --(23,4) ;
\draw [dashed]  (24,0) --(24,3) ;
\draw [dashed]  (25,0) --(25,2) ;
\draw [dashed]  (26,0) --(26,1) ;

\draw  (30,0)--  (36,2) --  (42,0)   ;
\draw  (36,2)--(36,6)  ;

\draw [dashed]  (32,0)--  (32, .66) --  (31,1)   ;
\draw  [dashed] (34,0) -- (34,1.32)--(32,2) ;
\draw [dashed] (36,0) --(36,2)--(33,3)  ;

\draw [dashed]  (36,2)--  (39,3)   ;
\draw  [dashed] (34,4) -- (36,3)--(38,4) ;
\draw [dashed] (35,5) --(36,4)--(37,5)  ;
\draw  [dashed] (38,0) -- (38,1.32)-- (40,2) ;
\draw [dashed] (40,0) --(40,.66)--(41,1)  ;

\end{tikzpicture}

\vskip 1cm \caption{Foliating a simplex}\label{simp}
\end{figure}


Again this map can be extended to every
subdivided $2$-simplex so that it commutes
with the $G$-action.
We change $X$ to be this subdivided complex.
Regard $X$ as a $2$-complex in which
each cell is attached via a loop of length 
three.  We can make a complex of groups
in which each $G_{\sigma }$ is the trivial
group.  Since $G$ is the fundamental group
of $X$ it is the fundamental group of this
complex of groups. We have described
a way of foliating the $2$-cells which
shows that this complex of groups resolves
$T$.

We now describe some moves on
a resolving $2$-complex which can be made
on a resolving complex which change a resolving 
$2$-complex to
another resolving $2$-complex.
\vfill
\noindent
{\bf Move 1}.  Subdividing a $1$-cell.

Let $\gamma $ be a $1$-cell, with vertices
$u, v$, which may be the same.  
This can be replaced by two $1$-cells $\gamma
_1, \gamma _2$ and a new vertex $w$,
so that $\gamma _1$ has vertices $u, w$ and
$\gamma _2$ has vertices $v, w$.
The groups associated with $w, \gamma _1,
\gamma _2$ in the new complex of groups are
all $G(\gamma )$.  The attaching maps of
$2$-cells are adjusted in the obvious way.

\vfill
\noindent
{\bf Move 2}. Folding the corner of a $2$-cell.

Suppose that one end of a foliated 
$2$-cell is as in Fig~\ref {fig:2cell}. Thus
$v$ is the end vertex of the $2$-cell
and adjacent vertices are $x, y$
and $x, y$ are mapped to the same point
of $T$, so that they lie on the same vertical
line. Let the adjacent $1$-cells to $v$
be $e_1$ and $e_2$, which conflicts with
our earlier notation but is in line with that
of \cite {[D2]} and \cite {[BF1]}.  Let the groups associated with
the cells (in the complex of groups) be
denoted by the corresponding  capital
letters.

Folding the corner results in a fold
of the graph of groups associated with
the $1$-skeleton of $X$.  Such a fold is
one of three types which are listed in
\cite {[BF1]} (as  Type A folds) or in \cite {[D2]}.  They
are shown in Fig~\ref{folds} for the reader's
convenience.  As the group acting is always $G$ it is not necessary
to carry out vertex morphisms
(see \cite {[D2]})  which are necessary when carrying out morphisms of trees rather than graphs.

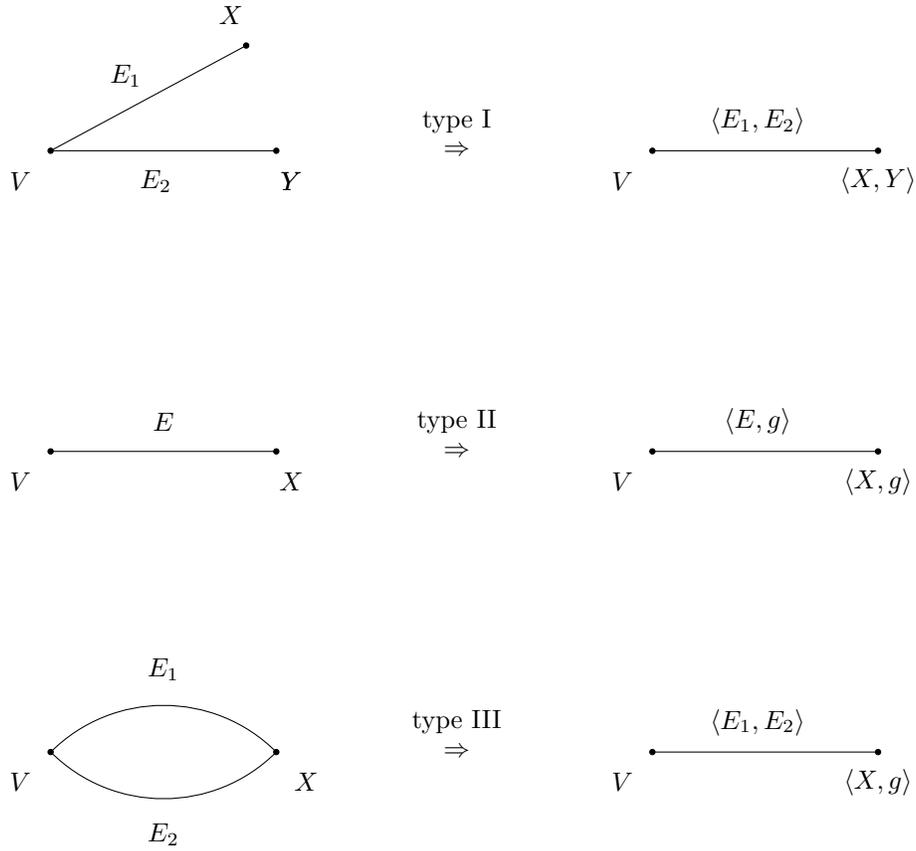
\begin{figure}

\centering
\begin{tikzpicture}[scale=.2]

     \filldraw (0,20) circle (5pt);
  \filldraw (15,20) circle (5pt);  
 
  \filldraw   (40,20) circle (5pt);
   \filldraw (55,20) circle (5pt);  
 
  \filldraw   (13,27) circle (5pt);
  \draw (0,20) -- (15,20) ;
\draw (0,20) --(13,27) ;
\draw (40,20) --(55,20) ;

\draw  (-2,18)  node {$V$} ;
\draw  (38,18)  node {$V$} ;
\draw  (5,25)  node {$E_1$} ;
\draw  (7,18)  node {$E_2$} ;
\draw  (12,29)  node {$X$} ;
\draw  (16,18)  node {$Y$} ;
\draw  (47,22)  node {$\langle E_1, E_2 \rangle$} ;
\draw  (55,18)  node {$\langle X, Y\rangle $} ;

\draw  (27,22)  node {type I} ;
\draw  (27,20)  node {$\Rightarrow $} ;
\filldraw (0,0) circle (5pt);
  \filldraw (15,0) circle (5pt);  
 
  \filldraw (40,0) circle (5pt);
    \filldraw (55,0) circle (5pt) ;

\draw (0,0) --(15,0) ;
\draw (40,0) -- (55,0) ;

\draw  (-2,-2)  node {$V$} ;
\draw  (38,-2)  node {$V$} ;

\draw  (7.5,2)  node {$E$} ;
\draw  (16,-2)  node {$X$} ;

\draw  (47,2)  node {$\langle E, g\rangle$} ;
\draw  (55,-2)  node {$\langle X, g\rangle $} ;
\draw (27,2) node {type II} ;
\draw (27,0) node {$\Rightarrow $} ;
\filldraw (40,-20) circle (5pt);
  \filldraw (55,-20) circle (5pt);  
 
  \filldraw (0,-20) circle (5pt);
    \filldraw (15,-20) circle (5pt) ;
    \draw (15,-20) arc (45: 135: 10.6) ;
   \draw (15,-20) arc (315: 225: 10.6) ;
\draw (40,-20)--(55,-20) ;

\draw  (-2,-22)  node {$V$} ;
\draw  (38,-22)  node {$V$} ;

\draw  (7.5,-14.5)  node {$E_1$} ;
\draw  (7.5,-25.5)  node {$E_2$} ;

\draw  (17,-22)  node {$X$} ;
\draw  (16,18)  node {$Y$} ;

\draw  (47,-18)  node {$\langle E_1,E_2\rangle$} ;
\draw  (55,-22)  node {$\langle X, g\rangle $} ;
\draw (27,-18) node {type III} ;
\draw (27,-20) node {$\Rightarrow $} ;

\end{tikzpicture}

\vskip 1cm \caption{Folding operations}\label{folds}
\end{figure}

The attaching word of the $2$-cell, whose
corner has been folded is changed 
in a way which we will describe in an example.
One can arrange that the joining element at the pivot vertex is trivial, by changing the lift of the spanning tree.
In this case, any other attaching word of a $2$-cell
that involves $e_1$ or $e_2$,
$\bar e_2$ is
replaced by the folded edge element $<e_1, e_2>$ and $e_1$
is replaced by $<e_1, e_2>$.
Let the new complex of groups be $X'$

Clearly there is a surjective homomophism 
$\phi : \pi (X, S) \rightarrow \pi (X', S')$
in which $g_ve_1g_v^{-1}$ and $e_2$
are both mapped to $<e_1, e_2>$. In fact this
homomorphism is an isomorphism 
since the resolving isomorphism 
$\alpha : \pi (X, S) \rightarrow G$ factors
through $\phi $.  We conclude that
$X'$ also resolves  the $G$-tree $T$.

If both the upper semi-circle and the lower semi-circle
consist of a single $1$-cell, then folding results in the elimination
of a $2$-cell, and a reduction in the number of $1$-cells.
 \

{\bf Move 3} Contracting a leaf.

Consider a foliated $2$-cell .  Let $\ell $ be a particular
vertical line of the foliation.  This will contain  points
$u, v$ of the upper semi-circle and lower semi-circle respectively.
After subdividing the relevant $1$-cells, it can be assumed that
these points are vertices.  Contracting the leaf $\ell $ results
in the $2$-cell $\sigma $ being replaced by two $2$-cells $\sigma
_1$ and $\sigma _2$.  The vertices $u, v$ become a single vertex $w$
and its group $G_w$ is the subgroup of $G$ generated by $G_u$ and
$G_v$ in $G$, except if $u,v $ belong to the same $G$-orbit, in which
case $G_w$ is generated by $G_u$ and an element $g \in G$ 
such that $gv = u$.  
Let $g_u, g_v$ be the respective elements of $G_u$ and $G_v$ in the attaching word for $\sigma$.
Let the edge after reaching $u$ in the attaching word end up in
$\sigma _2$.  This means that the edge after reaching $v$ ends up in $\sigma _1$. 
Suppose first that $u, v$ are in different orbits, then after the move the element for $w$ in $\sigma _2$ is $g_u$ and the element for $w$ in $\sigma _1$ is $g_v$.
Note that an edge has to be removed from the spanning tree $S$.
If $u, v$ are in the same orbit then the element for $w$ in $\sigma _2$
is $g_ug$ and the element for $w$ in $\sigma _1$ is $g_vg^{-1}$.
 
 A similar argument to that
for Move 2 shows that the complex we have created also resolves the $G$-tree $T$.
\

Let $\sigma $ be a $2$-cell of $X$.
We now examine what can happen as we 
repeatedly fold corners of $\sigma $,
at each stage replacing $\sigma $ by the
new $2$-cell created.
Since each $1$-cell of $\tilde X$ injects into 
$T$ we can assign each $1$-cell $\gamma $ of
$X$ a length,
namely the distance in $T$ between $\theta (u)$
and $\theta (v)$ where $u, v$ are the vertices
of a lift of $\gamma $ in $\tilde X$.

As above let $x$ be the corner vertex and
let $e_1, e_2$ be the incident edges.

 If $e_1,
e_2$ have the same length, then we can
fold the corner of $\sigma $.  If $e_1$
is shorter than
$e_2$ then subdivide $e_2$ so that the initial
part has the same length as $e_1$ and then
fold the corner.
If $e_2$ is shorter than $e_1$ then we
subdivide $e_1$ and then fold the corner.
Now repeat the process.
This process may terminate when
all the $2$-cell is folded away.

However it may happen that the folding
sequence is infinite i.e. it never terminates.

First we give an example making it easier to understand the following general explanation.
This example is a corrected  version of  Example 6 of \cite{[D]}.  
\begin {exam}\label {example}
\vskip 5mm
\begin{figure}[htbp]
\centering

\begin{tikzpicture}[scale=2]
          
    \path (0,5) coordinate (p1);
    \path (.707, 5) coordinate (p2);
    \path (1.707,5) coordinate (p3);
    \path (2.414,5) coordinate (p5);
     \filldraw (p2) circle (1pt);
  \filldraw (p3) circle (1pt);  
  \draw [left]  (0,5) node {$A$} ;
  \draw [below]  (p2) node {$B$} ;
  \draw [above]  (p3) node {$D$} ;
  \draw [right ]  (p5) node {$C$} ;
  \draw  (4, 5)--(4.353, 5.353)--(4.853, 5.353)--(5.353, 4.853)--(5,4.5)--(4.5,4.5)--cycle ;
  \draw  (4.353,4)--(4.753, 4.25)--(5.23, 4.25)--(5.353, 4.053)--(5,3.8)--(4.5,3.8)--cycle ;
    \draw  (4.5,3)--(4.753, 3.2)--(5, 3.3)--(5.23,3.203)--(5.353, 3)--(4.926, 2.8)--cycle ;
      \draw  (4.753, 2)--(5, 2.3)--(5.23,2.203)--(5.353, 2)-- (5.2, 1.85)--(4.926, 1.8)--cycle ;
      
      \draw [right]  (5.2,1.85) node {$_B$} ;
  \draw [below]  (4.65,2.1) node {$_B$} ;
 \draw [above]  (5,2.3) node {$_D$} ;
  \draw [above ]  (5.23,2.203) node {$_A$} ;
   \draw [right]  (5.353,2) node {$_D$} ;
  \draw [below ]  (4.926,1.8) node {$_C$} ;

\draw (.3, 5) node {$>$} ;
\draw [above] (.3, 5) node {$e$} ;
\draw (1.2, 5) node {$>$} ;
\draw [above] (1.2, 5) node {$f$} ;
\draw (2.1, 5) node {$>$} ;
\draw [above] (2.1, 5) node {$g$} ;
\draw (2.1, 4) node {$>$} ;
\draw (2.1, 3) node {$>$} ;
\draw (1.55, 3) node {$>$} ;
\draw [above] (1.55, 3) node {$h$} ;
\draw [above] (2.1, 3) node {$g$} ;
\draw [above] (2.1, 4) node {$g$} ;
\draw (1.1, 4) node {$<$} ;
\draw [above] (1.1,4) node {$e$} ;
 \draw [left]  (4,5) node {$_B$} ;
  \draw [below]  (4.353, 5.353) node {$_A$} ;
  \draw [below,left]  (4.5,4.5) node {$_D$} ;
  \draw [right ]  (5.353, 4.853) node {$_D$} ;
   \draw [above]  (4.835, 5.353) node {$_B$} ;
  \draw [right ]  (5, 4.5) node {$_C$} ;
  
  \draw [left]  (4.353, 4) node {$_A$} ;
  \draw [below]  (4.753,4.25) node {$_B$} ;
  \draw [above]  (5.23,4.25) node {$_A$} ;
  \draw [right ]  (5.353,4.053) node {$_D$} ;
   \draw [above]  (5,3.8) node {$_C$} ;
  \draw [below ]  (4.5,3.8) node {$_D$} ;
  
 \draw [left]  (4.5,3) node {$_D$} ;
  \draw [below]  (4.753,3.2) node {$_B$} ;
  \draw [above]  (5,3.3) node {$_D$} ;
  \draw [above ]  (5.23,3.203) node {$_A$} ;
   \draw [right]  (5.353,3) node {$_D$} ;
  \draw [below ]  (4.926,2.8) node {$_C$} ;

    \filldraw (p5) circle (1pt) ;
   \filldraw (p1) circle (1pt) ;

\draw (p1)--(p5) ;

\draw (1.614, 4) arc (0: 180: .2) ;
\draw (.907,5) arc (0: 180: .2) ;
\draw (1.907,3) arc (0: -90: .2) ;
\draw (2.3,2) arc (0: 180: .179) ;

    \path (0,4) coordinate (p1);
    \path (.707, 4) coordinate (p2);
    \path (1.707,4) coordinate (p3);
    \path (2.414,4) coordinate (p5);
    \draw [below]  (1.414,4) node {$A$} ;
  \draw [below]  (p2) node {$B$} ;
  \draw [above]  (p3) node {$D$} ;
  \draw [right ]  (p5) node {$C$} ;
     \filldraw (p2) circle (1pt);
  \filldraw (p3) circle (1pt);  
  \filldraw (1.414,4) circle (1pt);  

    \filldraw (p5) circle (1pt) ;
       

\draw (p2)--(p5) ;
   \path (1.414,3) coordinate (p1);
    \path (1.707, 2.707) coordinate (p2);
    \path (1.707,3) coordinate (p3);
     \path (2.12,2) coordinate (q3);
     \path (1.707,2) coordinate (q4);

    \path (1.414,2) coordinate (p4);
    \path (2.414,2) coordinate (p5);
   \path (2.414,3) coordinate (p6);
   
    \draw [below]  (p2) node {$B$} ;
  \draw [above]  (p3) node {$D$} ;
  \draw [right ]  (p5) node {$C$} ;
  
    \draw [below]  (q3) node {$B$} ;
  \draw [above]  (q4) node {$D$} ;
  \draw [right ]  (p6) node {$C$} ;

\draw [below]  (p1) node {$A$} ;
\draw [right ]  (p6) node {$C$} ;
     \filldraw (p2) circle (1pt);
  \filldraw (p3) circle (1pt);  
 
  \filldraw (q4) circle (1pt);
    \filldraw (q3) circle (1pt);  
 
  \filldraw (p4) circle (1pt);
\draw [above]  (p4) node {$A$} ;

    \filldraw (p5) circle (1pt) ;
       \filldraw (p6) circle (1pt) ;
   \filldraw (p1) circle (1pt) ;

\draw (p1)--(p6) ;
  \draw (p4)--(p5) ;    

\draw (1.707,2.707)--(p3);

\end{tikzpicture}
\vskip-3mm\caption{Folding sequence}\vskip-2mm
\end{figure}
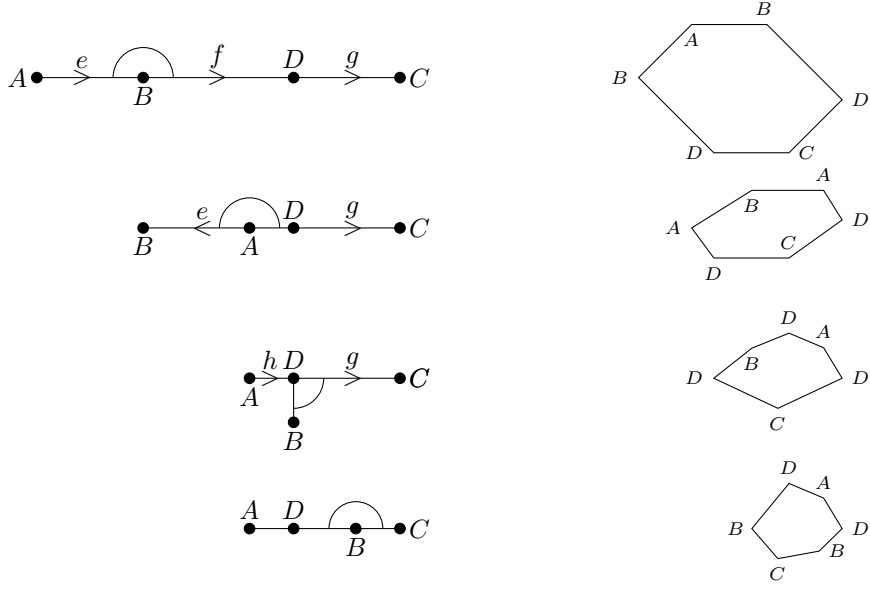

Let the complex $X_1$ have four  vertices $A, B, C, D$ and three oriented edges $e,f,g$.
Let $\iota e = A, \tau e = B, \iota f = B, \tau f = D, \iota g = D, \tau g = C$.
Let the groups of $A, B, C, D$ be finite cyclic of order $3$ and generated by $u, v, y, z$ respectively.
Let the $6$-sided $2$-cell be attached via the word  
$$w\cup w' =  \bar ea^{-1}ebfgc^{-1}\bar gd\bar f$$
Here $w = \bar ea^{-1}ebf$ and $w' = fd^{-1} gcg$.
In this case $X_1$ is a $2$-sphere with $4$ cone points.
Let $G_1$ be the group of this complex of groups.

The attaching word is describing a loop in $\tilde X_1$ the universal cover of
the complex of groups.  This loop maps to the loop,starting at $B$ $\bar eefg\bar g\bar f$ in $X_1$.
This loop is obtained by omitting the joining elements, which are elements of the
vertex group that has been reached at that point.  How a path in $\tilde X_1$ corresponds to 
such a word was described earlier.
We now discuss  how the  joining elements occur in $w\cup w'$.

In this case the $1$-skeleton of $X_1$ is a tree.  
We choose a particular lift of this tree in $X_1$ to the universal orbifold cover  $\tilde X_1$, which is  the (hyperbolic) plane tessellated by
$6$-gons.  Here $\tilde X_1$ is the universal cover of the complex of groups described above.
The attaching word traces out a  loop in $\tilde X_1$, which is the boundary of a fundamental region. Note that although the image of the path backtracks in $X_1$,
it is not allowed to backtrack in $\tilde X_1$.   This means that there must be non-trivial joining elements where the image backtracks.
In $\tilde X$ each fundamental region has $6$ vertices, including one point (incident with 3 edges in $\tilde X_1$) from the orbits corresponding to $A$ and $C$ and two vertices
(incident with $6$ edges in $\tilde X_1$)
from each of $B$ and $D$.
At one of the visits of the attaching word to vertices corresponding to $B$ (or $D$) we have to use a non-trivial  joining element.  We can choose where
this is.  
We get a presentation for $G_1$ in which the generators are $a, b, c, d$ and a relation
obtained by deleting the edges in the attaching word.  This is because the $1$-skeleton of $X_1$ is a tree.
Thus there is a relation $a^{-1}bd^{-1}c=1$.  There are also relations $a^3 = b^3 = c^3 = d^3 = 1.$

Clearly there is a surjective homomorphism $\phi $ from $G_1$ to the rectangle group  $$G = \langle a, b, c, d |a^3= b^3=c^3=d^3,  x = a^{-1}b =c^{-1}d, y = a^{-1}c = b^{-1}d
  \rangle .$$

Consider the folding sequence corresponding to the ``marking"   in which lengths are assigned to the edges  with $|e|  = |g| = 1,  |f|  = \sqrt 2$.
Initially we have the $2$-cell attached along $w\cup w' =  \bar ea^{-1}ebfgc^{-1}\bar gd\bar f$.
The attaching word is quadratic - its image in $X_1$ is $\bar eefg\bar g\bar f$ - and we will see that there is an infinite folding sequence
in which $w = \bar e a^{-1}ebf$  is folded against $\bar w'  = fd^{-1} gc\bar g$.
 The total length along top or bottom is $1 +\sqrt 2$.  After the first subdivision and fold we have a 
 new complex $X_2$ with the same vertices $A, B, C, D$ and with edges $b, c$ and a new  edge $h$ with length $\sqrt 2 -1$  with $\iota h = A, \tau h = D$  and the attaching word has become
 $eb\bar eb^{-1}a^{-1}bhgc^{-1}\bar gd\bar h.$   Note that the joining element $a^{-1}$ has changed to a conjugate $b^{-1}a^{-1}b$ as its position has changed.
 

 The $2$-cell has $w = eb\bar eb^{-1}a^{-1}bh, \bar w' = hd^{-1}\bar gcg$
After the next subdivision and  fold we have a new complex $X_3$ with the same vertex set but with  edges $h, g,j$ where $\iota j = D, \tau j = B$  and $j$ has  length $1 -(\sqrt 2 -1) = 2-\sqrt 2$  and the attaching word has become $jb\bar j(b^{-1}db)\bar h(b^{-1}a^{-1}b h  gc\bar g  $. Note that in this 
graph there is a vertex $D$ of valency $3$ whereas previously no vertex had valency more than $2$.
The attaching word visits the vertex $D$ three times.  As before we move the non-trivial  joining element so that it is not at the start or end point of $w = jb\bar j(b^{-1}db)\bar h (b^{-1}a^{-1}b)h$.  
As noted above this change of position of the  joining element corresponds to a change of the lift of a spanning tree  - in this case the whole of the $1$-skeleton $S_1$ of $X_1$.
Having chosen a lift $\tilde D$  of $D$ there are $27$ different lifts of $S_1$ to $\tilde X_1$.
These are acted on by the stabilizer of $\tilde D$ and there are nine different orbits under this action.
The attaching map must have at least one non-trivlal  joining element on a visit to $\tilde D$, since
otherwise one could have used the trivial group as the group at $D$.  We can choose the lift 
of $S_1$ so that the  joining element is non-trivial at exactly one visit.
In this case we do it so that the non-trivial  joining element is at a visit which is not the start or end point
of $w$.
We now have  $w = jb\bar j(b^{-1}db)\bar h (b^{-1}a^{-1}b)h$ 
and  $\bar w' = gc^{-1}\bar g$.  We can translate the whole lift by $b$ giving an attaching word
$jb\bar jd\bar ha^{-1}hg(bcb^{-1})\bar g$.   All we have done here is conjugate all the  elements by $b$ to make the  elements shorter.
The next subdivision and fold starts at $D$ and folds $j$ the shorter edge against $g$, so that we then have  a new edge $k$ with length $1 - (2-\sqrt 2) =\sqrt 2 -1$ replacing $g$. Here $\iota k = B, \tau k = C$
 and the attaching word is $\bar j d\bar h a^{-1}h j(c^{-1}bc) k(bcb^{-1}) \bar k$.   
 Now note that  the situation we have reached is similar to the initial
situation scaled by $\sqrt 2 -1$.    In fact the positions of $A$ and $C$ have been transposed from the original position.  To get an exact scaling carry out the next $3$ folds to get 
the initial position scaled by $(\sqrt 2 -1)^2$.

The foliation of $X_1$ corresponding to our marking, lifts to a foliation of $\tilde X$ and there is an $\R $-tree $T_1$ in which the points are leaves
of this foliation.  Clearly $T_1$ is a $G_1$-tree.
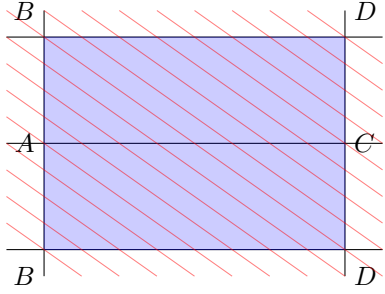
\begin{figure}[htbp]

\begin{tikzpicture}[xscale = .1, yscale = .0707, fill opacity =.25]
\centering

\draw (0,0)--(0,50) ;

\draw  (40, 0 )--  (40,50) ;

\draw (-5,5)--(45,5) ;
\draw (-5,25)--(45,25) ;
\draw (-5,45)--(45,45) ;

\filldraw [blue, opacity = .2]  (0,5)--(0, 45)--( 40,45) --(40,5) -- cycle  ;

\draw [fill opacity =1] [left] (0,0) node  {$B$} ;
\draw [fill opacity =1] [left] (0,25)  node  {$A$} ;
\draw [fill opacity =1] [left] (0,50) node  {$B$} ;
\draw [fill opacity =1] [right] (40,50) node  {$D$} ;
\draw [fill opacity =1]  [right] (40,25) node  {$C$} ;
\draw [fill opacity =1] [right] (40,0) node  {$D$} ;
  \draw [red, opacity = .5]  (5,0) --(-5,10) ;
\draw [red, opacity = .5]  (10,0) --(-5,15) ;
\draw [red, opacity = .5] (15,0) --(-5,20) ;
\draw [red, opacity = .5] (20,0) --(-5,25) ;

  \draw [red, opacity = .5]  (25,0) --(-5,30) ;
\draw [red, opacity = .5]  (30,0) --(-5,35) ;
\draw [red, opacity = .5] (35,0) --(-5,40) ;
\draw [red, opacity = .5] (40,0) --(-5,45) ;

  \draw [red, opacity = .5]  (45,0) --(-5,50) ;
\draw [red, opacity = .5]  (45,5) --(0,50) ;
\draw [red, opacity = .5] (45,10) --(5,50) ;
\draw [red, opacity = .5] (45,15) --(10,50) ;

  \draw [red, opacity = .5]  (45,20) --(15,50) ;
\draw [red, opacity = .5]  (45,25) --(20,50) ;
\draw [red, opacity = .5] (45,30) --(25,50) ;
\draw [red, opacity = .5] (45,35) --(30,50) ;

\draw [red, opacity = .5] (45,40) --(35,50) ;

\end{tikzpicture}
\caption{The image of a $2$-cell}\label{image}
\end{figure}
Let $C$ be the cube complex for the rectangle group  $G$.    As described in \cite{[D5]} there is a cube complex $\tilde C$ on which $G$ acts.
Thus $$G = \langle a, b, c, d |a^3= b^3=c^3=d^3,  x = a^{-1}b =c^{-1}d, y = a^{-1}c = b^{-1}d
  \rangle .$$
  There is a Euclidean subspace $E$  of $\tilde C$ acted on by $\langle x, y \rangle $
Assume that $x$ acts on $E$ by translation $2$ in the $x$-direction, and $y$ by translation of $2\sqrt 2$    There is a foliation on $C$ which induces a foliation on $\tilde C$ and $E$ given by the lines $x +y = c$.   The leaves of the foliation on $\tilde C$ give an $\R $-tree $T$.
  There is a map $\theta :  X_1 \rightarrow C$, which induces a  homomorphism,  denoted
$\theta  ':  G_1 \rightarrow G$.   The map $\theta $ induces a map $\tilde \theta : \tilde X_1 \rightarrow \tilde C$ in which a $2$-cell of $\tilde 
X_1$ maps into $E$ as indicated in Fig~\ref{image}.  It can be seen that the foliation of $E$ lifts to the foliation on  $\tilde X_1$ which is the one induced on the $2$-cell of $X_!$ corresponding to the marking with $|e| = |f| = 1, |g| = \sqrt 2$.
Thus there is a map $\bar \theta :  T_1 \rightarrow T$ which commutes with the actions of $G_1$ and $G$ via $\theta '$.

In $G_1$ the three  elements $y =b^{-1}d, y' = a^{-1}c$ and $ x = a^{-1}b = c^{-1}d$ freely generate a subgroup $F$ and the minimal $F$-subtree of $T_1$
is non-simplicial.   The elements $y, y'$ are hyperbolic elements with the same hyperbolic length $2\sqrt 2$.     Their axes intersect is a segment
of length $2\sqrt 2$.   In $G$ these elements become equal  and so have the same axis.   There may be a $G_1$-tree $T'$ for which
there are morphisms $T_1 \rightarrow T' \rightarrow T$, in which the two axes have a larger intersection than in $T_1$.

\end{exam}

\begin{proof} [Proof of Theorem 1.1]
Let $T$ be a $G$-tree, where $G$ is a finitely presented group.  We have seen that $G$ is the fundamental group of a complex of groups $G(X)$ that resolves the action.

Let $\sigma $ be a $2$-cell of $X$.
We now examine what can happen as we 
repeatedly fold corners of $\sigma $,
at each stage replacing $\sigma $ by the
new $2$-cell created.
Since each $1$-cell of $\tilde X$ injects into 
$T$ we can assign each $1$-cell $\gamma $ of
$X$ a length,
namely the distance in $T$ between $\theta (u)$
and $\theta (v)$ where $u, v$ are the vertices
of a lift of $\gamma $ in $\tilde X$.

As above let $x$ be the corner vertex and
let $e_1, e_2$ be the incident edges.

 If $e_1,
e_2$ have the same length, then we can
fold the corner of $\sigma $.  If $e_1$
is shorter than
$e_2$ then subdivide $e_2$ so that the initial
part has the same length as $e_1$ and then
fold the corner.
If $e_2$ is shorter than $e_1$ then we
subdivide $e_1$ and then fold the corner.
Now repeat the process.
This process may terminate when
all the $2$-cell is folded away.

However it may happen that the folding
sequence is infinite i.e. it never terminates.
We examine when this happens. 
Suppose this is the case and that 
the $2$-complexes in the sequence are 
$X_n , n = 1, 2, \dots \ \ $.

We can assign lengths
to the edges ($1$-cells)  of $X_n$.
Traversing the top semi-circular boundary of
the $2$-cell $\sigma $ determines
 to a path (or rather walk) $w$ in
the $1$-skeleton of $X_1$. Let
$w'$ be the path corresponding to the lower 
semi-circular boundary.
These paths
are usually not segments - they can even
backtrack. Let $\ell _n$ be the total length
of edges of
$X_n$.  It is clear that $\ell _n
\geq \ell _{n+1} \geq 0$.  We have $\ell _{n+1}
= \ell _n$ if and only if the fold is a
subdivision or a type II fold.
In going from $X_n$ to $X_{n+1}$
an arc $[y_n, y_{n+1}]$ of the upper
semicircular boundary of $\sigma $ is
identified with an are $[y_n', y_{n+1}']$
of the lower semicircular boundary.
Each such arc is identified with a $1$-cell
of $X_n$ and so has a length. 
In $X_n$ the folding has identified 
$[x, y_n]$ with $[x, y_n']$.
We assume that $y = \lim _n
y_n$, and that
$y' = \lim y_n'$.  It is possible that
$y = y'$ is the end point of $\sigma $ and
we will see that this is often the case.
Let 
$\lambda _n$ be the length of the arc
$[y_n, y]$. Thus $\lambda _n$ is the length
of the arc which remains to be folded.

We show that there can only be finitely
many type II folds in our sequence.
This is because there can only be
a finite number of type II folds to start
with as each such fold will use up the full
length of an edge of $X_1$. In our
sequence, a type I fold can only be followed
by a  type II fold if the type II fold is
between edges in the same orbits as as the
ones that  were folded together in the type I
fold. Thus there is a vertex in $[x,y]$ such
that the adjacent edges are in the same orbit. 
Such a vertex must have been a vertex in the 
original path $[x,y]$ in $X_1$ and so 
this happens only finitely many times.
Each type III fold  decreases the first 
Betti number of the quotient graph
and so there can only be a finite number of
type III folds. 
In our sequence there are therefore only 
finitely many type II or type III folds.
Assume then that all folds in the sequence 
are of type I. 
  Consider the 
subspace  of $X_1$ which is
the union of the images of the paths
corresponding to $[x,y]$ and $[x,y']$.  If
this is not a  subgraph of $X_1$, then
one of the paths corresponding to $[x,y],
[x,y']$ in $X_1$ must end in part of an
edge not visited by the other path. It is not
hard to see that this will not produce an
infinite folding sequence. Thus we assume that
this  subgraph is all of the $1$-skeleton
of $X_1$.

In our sequence of subdivision and type I
folds the number of edges in the quotient graph
does not increase, since any subdivision
which increases the number of edges by one
is immediately followed by a type I fold which
reduces it by one.  Clearly there can only
be a finite number of type I folds which
are not preceded by a subdivision, since the
number of such folds is bounded by the number
of edges of $X_1$.
It may happen that a fold at the 
$n$-th stage involves  an
edge  which is not in the subgraph
$X'_{n+1}$ determined by the remaining
folding sequence. This can happen for only a
finite number of folds, since if this happens
$X'_{n+1}$ has fewer edges than
$X_n$. Thus we assume that each folded
edge is in the subgraph determined by the
remaining folding sequence.

For a type I fold
$\ell _n$ and $\lambda _n$ are reduced 
by the same amount. 

Since we are assuming that
each folded edge is in the subgraph
determined by the remaining folding
sequence, it is 
clear that  $\ell _n$ tends to zero
as $\lambda _n$ tends to zero. 
Since $\ell _n - \lambda _n$ is constant,
it follows that $\ell _n = \lambda _n$.

  Let $w_y, w'_{y'}$ be the directed paths in
$X_1$ which are the images of
 $[x, y], [x, y']$
respectively. Clearly they are initial
parts of the paths
$w, w'$, so they  
 begin at the same
point.  In fact we can assume that they end at the
same point by using a Move 3 to contract the leaf that contains
the points $y$ and $y'$.   In fact we will show that $y$ and $y'$
are always vertices in the original graph.

From length considerations every edge of $X_1$ occurs exactly
twice in $w\cup w'$ or at least one
edge occurs only once.
If the latter occurs we will arrive at
a contradiction by showing that the folding
sequence 
must have been finite.
Let $e$ the edge which occurs only once in
$w\cup w'$.  Without loss of generality
suppose it is in $w$.  In fact we can 
assume that it is the first edge of $w$,
since we can fold away any edges which
precede it. This folding will not affect
the edge $e$.  There is also a folding sequence
starting at the other end of $\sigma $.  It
is not hard to see that this must also be 
an infinite sequence and in the limit
all of $w\cup w'$ is folded away.  Folding
away those edges which occur before
$e$ in this sequence and after $e$ in the
original sequence, we arrive at a new 
$2$-cell in which the entire path
$w$ consists of a single edge $e$. 
But such a folding sequence must be finite -
it will just fold $e$ onto the path
$w'$.
We have the desired contradiction.

An infinite folding sequence therefore
occurs when there is a $2$-cell in which
the attaching map contains every
edge exactly twice.  As we shall see, however,  a quadratic attaching map does not necessariy correspond to an infinite folding sequence.

We want to show that we can carry out folding on
the different $2$-cells and end up with a complex in which each cell
is attached via a quadratic word.

After carrying out a finite number of Type 3 moves we can assume that each leaf of the foliation
 intersects the top and bottom of each $2$-cell in at most one vertex.
 The argument above shows that  the limit points $y, y'$ of an infinite
 folding sequence must be vertex points on the same leaf of a foliation and so $y = y'$ will be an end point of the $2$-cell.  If one considers the folding sequence starting from the other end of the $2$-cell,
 we see that the first point reached where the attaching word becomes quadratic must also
 correspond to a leaf of the foliation which, if it was different from an end point of the $2$-cell, would contain two vertices. Thus every $2$-cell corresponds to a finite folding sequence or it corresponds to
 an infinite folding sequence given by a quadratic attaching word.
 
Suppose  a complex is given by a single quadratic word $w \cup w'$.  A {\it marking } is an assignment
of positive lengths to the letters in such a way that the total length of $w$ is the same as that of $w'$.
Let $a_1, a_2, \dots , a_r$ be the letters which lie both in $w,$ and $w'$.   Let
$b_1, b_2, \dots b_s$ be the letters which occur twice in $w$ (and so not in $w'$) and let $c_1, c_2, \dots , c_t$
be the letter which occur twice in $w'$, then the $a_i b_j, c_k$ can be assigned arbitrary positive lengths
$\alpha _i, \beta _j, \gamma _k$ subject only to the single constraint $\beta _1 +\beta _2 \dots  + \beta _s
= \gamma _1 + \gamma _2 \dots + \gamma _t$.  The subspace of the real numbers generated by the coefficients therefore, has maximal dimension $r $ if there are no letters that occur twice in either  $V$ or $W$ and it has dimension $r + s + t -1$ if there are letters that do occur twice in either $w$ or $w'$.

Choose a resolving complex $X$ that has fewest $1$-cells.  Each attachment of a $2$-cell must induce an infinite folding sequence,  since Type I and Type III folds result in a reduction in the number of edges,
so any Type I fold must be preceded by a subdivision, and there are no Type III folds.
We define an equivalence relation on the set  $U$ of $1$-cells that occur as a face  of a $2$-cell of $X$.   We require that $e \sim f$ if there is a $2$-cell that includes both $e, f$ in its attaching map.  We take $\sim $ to be the smallest equivalence relation for which this is the case.
For any $2$-cell $\sigma $ of $X$   all the $1$-cells to which it is attached lie in a single equivalence class.   Thus for each equivalence class there
is a subcomplex consisting of those $1$-cells and its vertices together with the $2$-cells attached to that class.   Any two such complexes intersect  in
a set of vertices, but no edges.    
If two $2$-cells $\sigma , \sigma '$ share an edge $e$  and vertex   $u$ in their  boundaries, then we will see that we can choose the same  joining elements in $G_u$ for the 
two $2$-cells.
This is because,
as in the example above, if a vertex requires a joining element, then
at some  stage in the folding sequence the attaching words for both $\sigma$ and $\sigma '$  will contain a subword  of the form $fj_u\bar f$.     In the tree $T$,  $f$ will map to an arc,
$\tau f = u $ will map to a point  $v$ and $j_uf $ will map to an arc intersecting $f$ in the single point $v$.    Thus $f $ and $j_uf$ determine different {\it directions }  $d_1, d_2$  at $v$, such that $d_2 = j_vd_1$.    This will be true for both the attached two cells, so that we can choose the same $j_u$ for both  attaching words.

For the moment let us assume that there is a single equivalence class, and so there is a single subcomplex $X$ itself.

In the resolution of the action of $G$ on $T$ each $1$-cell $e$ is effectively assigned a length $|e| \in \R$, which is the length of the arc
in $T$ joining the images in $T$ of the end points of a lift of $e$ to the universal cover $\tilde X$.
Let $A$ be the the subgroup of $\R $, regarded as an additive group, generated by the set $\{ |e| |  e \in X^1 \} $.
The group $A$ is isomorphic to $\Z ^n$ for some $n$.     
Let $P_n$ be a parallelepiped group corresponding to an $n$-cube, in which we will assign orders to the vertex elements in a certain way.

As described in \cite {[D5]}  the group  $P_n$ acts on a $1$-connected, $n$-dimensional cubing $C_n$ that contains an $n$-dimensional Euclidean space
$E_n$ and $P_n$ contains a free abelian rank $n$ subgroup $J_n$ that acts on $E_n$ by translations of $2$ units in each of the coordinate directions.
The space $P_n\backslash C_n$ is obtained from $J_n\backslash E_n$ be identifying a single  antipodal pair of $n$-cells in $J_n\backslash E_n$.
We show that there is a subgroup  $G'$ of $G$ generated by cyclic subgroups  of  distinct vertex groups of $X$  and a map $\theta : \tilde X \rightarrow  C_n$ which is equivariant with respect to a homomorphism $G' \rightarrow P_n$.
For each $2$-cell $\sigma $ in $X$ there is a lift $\tilde \sigma $ such that $\theta (\tilde \sigma ) \subset E_n$, and the map $\theta $ is defined by specifying how $\theta $ acts on these $2$-cells.     

Suppose $A$ is generated by the real numbers $\alpha _1, \alpha _2, \dots , \alpha _n$.   We assume now  that $J_n$ acts on $E_n $ by translations of $2\alpha _i$ in each of the coordinate directions.    We give $E_n$ the structure of a cell complex in the obvious way so that, as for the rectangle group, 
$J_n$ acts cellularly and there is one orbit of  $n$-cells subdivided into $2^n$ orbits of smaller $n$-cells.
For the left hand  vertex $v$ in a particular $2$-cell $\tilde \sigma $ of $\tilde X$ let
$\theta (v)$ be the origin in $E_n$.   Each $1$-cell in $\tilde X$ has a particular length in $A$, and this length will determine a vertex of
$E_n$.    Proceeding around the boundary of $\tilde \sigma $ will determine a loop in the positive quadrant of $E_n$.  The distance from
the origin will increase as one passes along the top or bottom of $\tilde \sigma $ away from $v$ and one will reach the same point which is
the image of the right hand vertex of the $2$-cell.  If a vertex is visited more than once in passing along the top or bottom, then on one of the visits one has to use a joining element to pass to the antipodal subcube in $E_n$.  This will mean that the path traced out in $E_n$ never backtracks, though its image in $X$ will backtrack.

We illustrate the above argument with another  example.
\begin {exam} \label {example2}
Consider the pair $(a\bar ab\bar bc\bar c,  d\bar de\bar e)$.   Suppose a $2$-cell  corresponding to this pair arises in the action of a group $G$ on an $\R $-tree with marking 
$|a|= 1, |b| = \sqrt (2). |c|= \sqrt (3), |d|=\sqrt (5). |e| = 1 +\sqrt (2)+\sqrt (3) -\sqrt (5)$.

In this case $n =4$ .  The path traced out by the word $ww'$ visits the vertices $(0,0,0,0),(1,0,0,0), (2,0,0,0), (2,1,0,0), (2,2,0,0), (2,2,1,0), (2,2,2,0), (1,1,1,,1), $

$(0,0,0, 2), (0,0,0,1), (0,0,0,0)$.

\end {exam}

If the $1$-skeleton of $X$ is a tree  $X^1$,  in the loop in $X^1$ corresponding to the boundary of $\sigma $ the path corresponding to 
successive visits to a particular vertex will pass over each edge an even number of times.   It the edge is oriented then it must pass over the
edge the same number of times in each direction.   This means that in $E_n$ if two vertices of $\tilde \sigma $ are in the same $G$-orbit,
then their images in $E_n$ are in the same $J_n$-orbit.   Consider the subgroup $G_{\sigma}$ of $G$ generated by the joining elements of $\sigma $.
The group $G_{\sigma} $ is the  fundamental group of the complex  $X_{\sigma}$  of groups corresponding to $\sigma $.   All the vertex groups and edge groups 
are cyclic.  Each one is generated by a power of a joining element.   An element that fixes an edge of $\tilde X_{\sigma }$ must fix every edge for the 
reason explained above.
Thus $G_{\sigma }$ has a cyclic normal subgroup $N_{\sigma}$ such that $G_{\sigma}/N_{\sigma}$ acts on $\tilde X_{\sigma }$ with trivial edge stabilizers.   We now show that there is a  a map $\theta _{\sigma} : \tilde X_{\sigma} \rightarrow \tilde C_n$ which is equivariant with respect to a homomorphism from $G_{\sigma }$ to $P_n$ with kernel $N_{\sigma }$. 

   We map $G_{\sigma }$ into $P_n$ by mapping $N_{\sigma}$ to the identity element and giving each joining element to a vertex element in which
   its order is the order of that element modulo $N_{\sigma }$.
The defining  relations between the joining elements of $G_{\sigma }$ are given
by the attaching maps of $\sigma $ as described above.    If we map $\tilde \sigma $ into $E_n$ then the relation is a consequence of
the relations of $J_n$. Thus we have a homomorphism from $G_{\sigma }$ to $P_n$.

If a different $2$-cell $\sigma '$ of $X$ shares an edge $e$ with $\sigma $, then there will be a lift $\tilde \sigma '$ that shares an edge with
$\tilde \sigma$ and the boundary map of $\sigma '$ will determine a closed path in $E_n$.   Thus the maps $\theta _{\sigma}$ and $\theta _{\sigma '}$  match up nicely and carrying out the extension to every $2$-cell we see that there will be a map $\theta : \tilde X \rightarrow C_n$ which
restricts to $\theta _{\sigma }$ on each $2$-cell $\tilde \sigma $.  This map will be equivariant with respect to $G'$, the subgroup of $G$ generated by 
all the $G_{\sigma}$  for every $2$-cell $\sigma $.
An infinite folding sequence will produce an in infinite non-decreasing sequence of edge groups whose union will be a normal subgroup of the group of the cube complex that is the kernel of the map to the target group.

 Note that any folding sequence results in a sequence of complexes that resolve the action
on $T$ and it can never be the case that the joining element becomes trivial in the folding.    
This is because at some stage in a folding sequence the joining element will lie between an edge $e$ and $\bar e$ and if the joining element
is trivial, then the action on $T$ will not be resolved.
Two vertices in different $G$ orbits in $\tilde X$
may end up in the same $J_n$-orbit (I don't know if this can happen - it may be that if two vertices are mapped to the same $J_n$-orbit, then some
folding sequence will result in  a Type III fold and the images of the vertices lying in the same $G$-orbit ).   If it can happen, then a way of dealing with this  is to give the vertex element in
$P_n$ as its order the lowest common multiple of the finite  orders of any  joining elements  mapped to it (modulo the smallest power of that element that fixes an edge) and map each joining element
to an appropriate power of the vertex element in $P_n$.

If the one skeleton $X^1$  of $X$ is not a tree, then let $W$ be a spanning tree  for $X^1$.   In this case we take $G'$ to be the subgroup of $G$ generated
by  the joining elements corresponding to a lift of $W$ to $\tilde X$ together with a connecting element generator for each edge
$e$ of $X$ that is not in $W$.  If $u, u'$ are the vertices of $e$ then the lift of  $W$ to $\tilde X$ will contain unique lifts  $\tilde u, \tilde u'$ of $u, u'$.
There will not usually be an edge of $\tilde X$ joining $\tilde u, \tilde u'$ but there is a lift  $\tilde e$ of $e$ with $\iota \tilde e$ = $\tilde u$.
The generator corresponding to $e$ is an element $c(e)$ of $G$ such that $c(e)^{-1}\tau \tilde e = \tilde u'$.    The edge $e$ is given a length ${1\over 2} |e|$ in our action
on $T$.
We want the corresponding generator to be mapped to a translation by $|e|$  in $J_n$.   In this case let 
$A$ be the subgroup of $\R $ generated by $\{ |e| | e\in EW \} \cup \{ {1\over 2} |e| | e \in EX^1 \setminus EW\}$.    Taking a generator $ {1\over 2} |e|$ for $ e \in EX^1 \setminus EW $  means that  there is a translation of  $|e|$ in $J_n$ since it is through an even number of units.
Let $G'$ be the subgroup of $G$ generated by all the joining elements and connecting elements.
We can now    define a homomorphism from $G'$ to $P_n$.   We map each joining element as before, as it will correspond to a vertex of the spanning tree.   For each edge $e$ that is not in the spanning tree,  we introduce a vertex that is the midpoint of the subdivided edge.   In $P_n$
we give the corresponding vertex element  $v(e)$ the order two.   If the vertex element in $P_n$ corresponding to the initial vertex of $e$ in $Y$ is $u$, then we map the connecting element $c(e)$ to the element $v(e)u^{-1}$,  which will then correspond to a translation of length $|e|$.

Now consider the case when $U$ may have more than one 
equivalence class for the relation $\sim $.   Let $Y$ be the graph in which $VY$ is the union $VX$ with  the  set of equivalence classes  $U / \sim $ that have more than one element. 

 Let
$(\Y, Y)$ be the graph of groups in which for each $v  \in VX$, $\Y (v) = G_{v'} $, where $G_{v'}$ is the $G$-stabiliser of the image $v'$ under $\theta$ of the lift of $v$ in the lift of $W$ to $\tilde X$. For each $v$ that is an equivalence class  in $U /\sim$ we take $\Y (v)$ to be the group $G'$ defined above.   The set $EY$ is the union of $EX \setminus U$  with 
an edge   for any $v \in VX$ that is the vertex of an edge $e$ that lies in the equivalence class $[e]$
of $U$, joining $v$ to $[e]$.  The group in $(\Y, Y)$ attached to this edge will be the cyclic subgroup of $G$ generated by the joining element of $v$.  If $v$ has not had a joinig element attached to it, then let the group attached to this edge be the identity subgroup.

The graph of groups we have constructed has the properties listed in the statement of the theorem, and so the proof is complete.
\end{proof}

Not every quadratic word will correspond to an infinite folding sequence.
We say that the pair $(V, W)$ of words  is {\it admissible} if  the word $V \cup \bar W$ is quadratic and 
for some marking the corresponding folding sequence starts at one end and finishes at the other.
Here we may have to include  joining elements to represent the attaching word in $\tilde X$, for example  if either $V$ or $W$ contains a subword $e\bar e$.
If $(V, W )$ is admissible then any marking which produces such a folding sequence is called
an {\it admissible} marking.
We  now explore which pairs of  words $(V, W)$ are admissible.
 Which pairs are admissible seems quite tricky to determine. 
 We have seen in  Example \ref {example} that  the pair $(\bar eef, fg\bar g)$ is admissible with the admissible marking $|e|= |g| =1, |f| =\sqrt (2$)
 On the other hand the pair $(aabbcc, ddee)$ of  Example \ref {example2} with marking $|a| =1, |b|=\sqrt (2), |c| = \sqrt (3), |d| = \sqrt (5), |e| = 1+\sqrt (2) +\sqrt(3) - \sqrt (5) $ reaches a Type III fold after eleven folds and the pair  $(aabbcc, ddeeff),  |a| = 1, |b|= \sqrt (7), |c|= \sqrt (11), |d| = \sqrt (2), |e| = \sqrt (3),  |f| = 1 +\sqrt (7) +\sqrt (11) - \sqrt(2) -\sqrt (11) $
 reaches a Type III fold at the hundredth fold.  That I am able to determine that the latter markngs are  not admissible  is thanks to Andrew Bartholomew for producing a programme that carries out folding sequences.  In the last two cases, there are  homomorphisms to discrete subgroups $\R ^4$  and $\R ^5$.
 
 I think the following is true.

\begin {conj} \label {admiss}  If $(V, \bar W)$ is admissible, then any marking of maximal dimension  is an admissibe marking.
\end {conj}

We have seen above that the pair $(a\bar ab,bc\bar c)$ is admissible.  Another example is $(abcd, dcba)$.  If the conjecture is true then
 $(aabb, cc)$,   $(aabb, ccdd)$. and $(ab\bar a\bar b, cd\bar c\bar d)$ are not admissible, as there are markings of maximal dimension that result in folding sequences that 
 give Tyoe III folds.

\begin {prop}  If Conjecture \ref {admiss} is true then so are the following statements.
\begin {itemize} 
\item [(i)] A quadratic pair is admissible if and only if one marking  of maximal dimension is admissible.
\item [(ii)] If a pair $(V, W)$ is admissible then any pair of edges which occurs in both $V$ and $W$ 
 must have the same orientation.
 \item [(iii)] Let $(V, W)$ be an admissible pair. If an edge pair with the same orientation occurs in $V$, then no edge pair with opposite orientations can occur in $W$.
 
\end{itemize}

 \end{prop}

\begin{proof}

 [(i)].  This is immediate.

[(ii)] Suppose the pair $(UaV, W\bar aX)$ is admissible.
By folding from both ends we can assume that $U$ and $X$ are  empty.   Let $b$ be the first edge of $W$, so that
$W = bW'$.
if $a$ is longer than $b$ then folding will give $(a'V, W'\bar a'\bar b)$ which means that
$(aV, W'\bar a\bar b)$ is admissible.   There will certainly be a marking of maximal dimension in which $a$ is longer than $b$.   Note that $W'$ has fewer edges than $W$.  By repeating this process
we get an admissible pair $(aV, \bar aU)$.  But since the first internal vertices match up, this cannot
be admissible.

[(iii)]  First note that if $V = aUaX$ and $W = bY\bar bZ$, then by folding and assuming $b$ is 
longer than $a$ produces an admissible pair $(UaX, b'Y\bar b'\bar aZ)$, which contradicts (ii).  Thus the original pair was not admissible.
In general if $W$ contains a pair of edges $b, \bar b$  that do not have the same orientation, then  by folding from one end we can assume that $V = aUaX, W = W'bY\bar bZ$. Then by assuming
the length of $a$ is greater than the total length of $W'$, by folding we get 
an admissible pair $(a'UW'a'X, bY\bar bZ)$.  But we have just seen that this cannot be admissible.
\end{proof}

\end{section}

\end{document}